\documentclass{amsart}

\theoremstyle{plain}
\newtheorem{thm}{\indent\sc Theorem}[section]
\newtheorem{prop}[thm]{\indent\sc Proposition}

\newtheorem{cor}[thm]{\indent\sc Corollary}
\theoremstyle{definition}
\newtheorem{defn}[thm]{\indent\sc Definition}

\theoremstyle{remark}

\def\proofwp{{\sc Proof}}

\numberwithin{equation}{section}
\def\cA{\mathcal{A}}

\def\cC{\mathcal{C}}

\usepackage{amscd}
\usepackage{amsmath,amssymb}
\usepackage{amssymb}
\usepackage{amsmath}
\usepackage{amsthm}
\usepackage{latexsym}
\usepackage{amsfonts}
\usepackage{color}
\usepackage{graphics}
\usepackage[all,knot]{xy}
\usepackage[all]{xy}
\xyoption{arc}
\usepackage{enumerate}
\usepackage{amstext}

\def\St{\mbox{\bf St }}

\newcommand{\RR}{\mathbb{R}}

\begin{document}

\title[Sheaf of smooth forms on transitive Lie algebroids]{On the sheaf of smooth forms on Lie algebroids over triangulated spaces}



\author{{\textsc{Jose R. Oliveira}}}
\address{Department of Mathematics, Minho University,
    Braga, Portugal}
\curraddr{}
\email{jmo@math.uminho.pt}
\thanks{The authors are partially supported by MICINN, Grant MTM2014-56950-P}



\date{}

\dedicatory{}


\maketitle






\begin{center}

\vspace{3mm}

\textsc{Abstract}

\vspace{3mm}

\end{center}


Cohomology of a compatible family of Lie algebroids defined on a family of transverse manifolds is defined. A sheaf of differential forms on a compatible family of Lie algebroids defined over regular open subsets of a simplicial complex is constructed. It is proved that sheaf is fine.


\vspace{3mm}







\section{Introduction}

\vspace{3mm}

Mishchenko and Oliveira in their paper \cite{mish-oli} considered families of transitive Lie algebroids obtained by restriction of a Lie algebroid to the simplices of the base, in which the base is smoothly triangulated by a simplicial complex. They defined a cochain algebra of differential forms in similar way to the ones presented in \cite{suli-inf} by Sullivan or in \cite{wity-git} by Whitney and proved that the cohomology of such algebra is isomorphic to the cohomology of the Lie algebroid considered. 

The work developed by Mishchenko and Oliveira led the author of the present paper to consider in \cite{jose-MV} compatible families of transitive Lie algebroids defined not only on triangulated manifolds but on general simplicial complexes.

In the present paper, we consider compatible families of transitive Lie algebroids defined over a family of manifolds with transverse intersections in an ambient space. The construction of the cochain algebra of piecewise forms follows the construction developed in \cite{mish-oli} or in \cite{jose-MV}. We focus our work mainly on compatible families of transitive Lie algebroids which are defined over the family made of all regular open subsets of a simplicial complexe. These families of Lie algebroids allow us to construct a sheaf of Lie algebroids. We will prove that the correspondent sheaf is fine.

\vspace{3mm}

Throughout the paper, all manifolds are smooth, finite-dimensional and possibly with boundaries of different indices. All simplicial complexes considered are finite and geometric. Simplex means always closed simplex. For each simplicial complex $K$, its geometric realization will be denoted by $|K|$.

If $M$ is a smooth manifold, $TM$ the tangent bundle to $M$ and $\Gamma(TM)$ the Lie algebra of the vector fields on $M$, a Lie algebroid on $M$ is a vector bundle $\pi:\cA\longrightarrow M$ with base $M$ equipped with a vector bundle morphism $\rho:\cA\longrightarrow TM$, called anchor of $\cA$, and a structure of real Lie algebra on the vector space $\Gamma(\cA)$ of the sections of $\cA$ such that the map $\rho_{\Gamma}:\Gamma(\cA)\longrightarrow\Gamma(TM)$, induced by $\rho$, is a Lie algebra homomorphism and the action of the algebra $\cC^{\infty}(M)$ on $\Gamma(\cA)$ satisfies the natural condition: $$[\xi,f\eta]=f[\xi,\eta] + (\rho_{\Gamma}(\xi)\cdot f)\eta$$ for each $\xi$, $\eta$ $\in \Gamma(\cA)$ and $f\in \cC^{\infty}(M)$. The Lie algebroid $\cA$ is called transitive if the anchor $\gamma$ is surjective. Let $\varphi:N\hookrightarrow M$ be a submanifold, possibly with boundaries of different indices and assume that $\cA$ is transitive. We recall that the Lie algebroid restriction of $\cA$ to the submanifold $N$, denoted by $\cA^{!!}_{N}$, is the Lie algebroid $\varphi^{!!}\cA$ constructed as inverse image of $\cA$ by the mapping $\varphi$ (see \cite{kuki-chw}, \cite{makz-lga} and \cite{mish-oli} for more details). The papers \cite{mish-cla}, \cite{mish-latb} and \cite{mish-obs} contain a summary of those issues.

\vspace{6mm}

\textbf{Acknowledgments}. I want to thank to Aleksandr Mishchenko, James Stasheff, Jesus Alvarez and Nicolae Teleman for their strong dynamism to discuss several topics concerning cohomology of Lie algebroids.


\vspace{6mm}

\begin{center}
\section{Algebra of piecewise smooth forms}
\end{center}

\vspace{3mm}

Let $K$ be a simplicial complex. A \textit{complex of Lie algebroids} on $K$ is a family $\underline{\cA}=\{\cA_{\Delta}\}_{\Delta\in K}$ such that, for each $\Delta\in K$, $\cA_{\Delta}$ is a transitive Lie algebroid on $\Delta$ and, if $\Delta$ and $\Delta'$ $\in K$ are two simplices of $K$, with $\Delta'$ face of
$\Delta$, the Lie algebroid restriction of $\cA_{\Delta}$ to $\Delta'$ is the Lie algebroid $\cA_{\Delta'}$. Following the idea of differential forms on cell spaces given in Whitney book's \cite{wity-git} or in Sullivan's work \cite{suli-inf}, a piecewise smooth form on the complex of Lie algebroids $\underline{\cA}=\{\cA_{\Delta}\}_{\Delta\in K}$ is a family $\omega=(\omega_{\Delta})_{\Delta \in K}$ such that, for each $\Delta\in K$, $\omega_{\Delta}\in \Omega^{p}(\cA_{\Delta};\Delta)$ is a smooth form on $\cA_{\Delta}$ and, if $\Delta'$ is a face of $\Delta$, $(\omega_{\Delta})_{/\Delta'}=\omega_{\Delta'}$ (cf. with \cite{mish-oli} or \cite{jose-MV}).

The direct sum $$\Omega^{\ast}(\underline{\cA};K)=\bigoplus_{p\geq 0}\Omega^{p}(\underline{\cA};K)$$ of all piecewise smooth forms on $\underline{\cA}$, equipped with the exterior product and the exterior derivative by the corresponding exterior product and exterior derivative on each algebra $\Omega^{\ast}_{ps}(\cA_{\Delta};\Delta)=\bigoplus_{p\geq 0}\Omega^{p}(\cA_{\Delta};\Delta)$, is a cochain algebra defined on $\RR$.

Let $L$ be a simplicial subcomplex of $K$ and $\underline{\cA}^{L}=\{\cA_{\Delta}\}_{\Delta\in L}$ the complex of Lie algebroids given by restriction of $\underline{\cA}$ to $L$ (see \cite{jose-MV} for definition of complex of Lie algebroids restriction). If $\omega=(\omega_{\Delta})_{\Delta\in K}\in\Omega^{p}(\underline{\cA};K)$ is a piecewise smooth form, we can define the restriction of $\omega$ to the subcomplex $L$, denoted by $\omega_{/L}$, to be the form $$\omega_{/L}=(\omega_{\Delta})_{\Delta\in L}\in\Omega^{p}(\underline{\cA}^{L};L)$$ The equality $d(\omega_{/L})=(d\omega)_{/L}$ holds. We obtain a new cochain complex of forms, the cochain complex $\Omega^{\ast}(\underline{\cA}^{L};L)$. For each $p\geq 0$, denote by $$r^{p^{K}}_{_{L}}:\Omega^{p}\big(\underline{\cA};K\big)\longrightarrow \Omega^{p}\big(\underline{\cA}^{L};L\big)$$ the map induced by restriction, that is, for each $\omega\in \Omega^{p}(\underline{\cA};K)$,  $$r^{p^{K}}_{_{L}}(\omega)=\omega_{/L}$$

\vspace{3mm}

\begin{prop} Keeping these hypotheses and notations as above, the following properties hold.

\begin{itemize}
\item For $p=0$, $r^{0^{K}}_{_{L}}:C_{ps}(|K|;\RR)\longrightarrow C_{ps}(|L|;\RR)$ is a homomorphism of algebras, in which $C_{ps}(|K|;\RR)$ denotes the
algebra over $\RR$ made by all continuous maps $f:|K|\longrightarrow\RR$ that are compatible with restrictions to the faces of $K$ and with smooth restrictions to the faces of $K$.

\item For each $p\geq 0$,  $r^{p^{K}}_{_{K}}=id_{\Omega^{p}(\underline{\cA};K)}$.

\item The map $r^{\ast^{K}}_{_{L}}:\Omega^{\ast}(\underline{\cA};K)\longrightarrow \Omega^{p}(\underline{\cA}^{L};L)$ is a morphism of cochain algebras.

\item If $T$ is a simplicial subcomplex of $L$ and $\underline{\cA}^{T}=(\cA_{\alpha})_{\alpha\in T}$ then, for each $p\geq 0$, $r^{p^{K}}_{_{T}}=r^{p^{L}}_{_{T}}\circ r^{p^{K}}_{_{L}}$ and so the diagram

$$
\xymatrix{
\Omega^{\ast}(\underline{\cA};K)\ar[dr]_{r^{\ast^{K}}_{_{T}}}\ar[rr]^{r^{\ast^{K}}_{_{L}}} & & \Omega^{p}(\underline{\cA}^{L};L)\ar[dl]^{r^{\ast^{L}}_{_{T}}} \\
& \Omega^{p}(\underline{\cA}^{T};T) &
}
$$ is a commutative diagram of cochain algebras.
\end{itemize}
\end{prop}

\vspace{3mm}

By the proposition 3.5 of \cite{jose-MV}, the map $r^{p^{K}}_{_{L}}$ is surjective.

\vspace{6mm}

\section{Sheaf of piecewise smooth forms}

\vspace{3mm}

In next section, we are going to consider a generalization of the concept of piecewise smooth cohomology given in the previous section by taking any family of smooth manifolds with transverses intersections. This will allow to construct a sheaf of the piecewise smooth forms which will be fine. As remarked in the introduction, all simplicial complexes considered are finite and geometric. Simplex means always closed simplex. For each simplicial complex $K$, its geometric realization will be denoted by $|K|$.

\vspace{3mm}

\begin{defn}\label{transfam} Let $\underline{K}=\{N_{1},\dots,N_{s}\}$ be a finite family of submanifolds of a smooth manifold $M$. The family $\underline{K}$ is said to be transverse if all intersections $N_{j_{1}}\cap \dots \cap N_{j_{e}}$, for any $j_{1}$, $\dots$, $j_{e}$ $\in \{1,\dots,s\}$, are transverse in the ambient manifold $M$ (cf. with \cite{gang-diff}).
\end{defn}


\begin{defn} Let $\underline{K}=\{N_{1},\dots,N_{s}\}$ be a transverse family of submanifolds of a smooth manifold $M$. A complex of Lie algebroids on $\underline{K}$ is a family $\underline{\cA}=\{\cA_{j}\}_{j\in \{1,\dots,s\}}$ such that, for each $j\in \{1,\dots,s\}$, $\cA_{j}$ is a transitive Lie algebroid on $N_{j}$ and, for each $i,j\in \{1,\dots,s\}$, one has $$(\cA_{j})^{!!}_{N_{j}\cap N_{i}}=(\cA_{i})^{!!}_{N_{j}\cap N_{i}}$$
\end{defn}

It is obvious that, by transitivity of restrictions of Lie algebroids, for any subset $\widetilde{J}$ of $\{1,\dots,s\}$ and any partition $\{\{j_{1},\dots,j_{r}\},\{i_{1},\dots,i_{t}\}\}$ of $\widetilde{J}$, we have $$(\cA^{!!}_{N_{r}})^{!!}_{N_{t}}=(\cA^{!!}_{N_{t}})^{!!}_{N_{r}}$$ in which $N_{r}=N_{j_{1}}\cap \dots \cap N_{j_{r}}$ and $N_{t}=N_{i_{1}}\cap \dots \cap N_{i_{t}}$.

\vspace{3mm}

Keeping the same hypotheses and notations from previous definition, we give now the definition of piecewise smooth form on the complex of Lie algebroids $\underline{\cA}$.

\vspace{3mm}

\begin{defn} A piecewise smooth form of degree $p$ ($p\geq0$) on $\underline{\cA}$ is a family $\omega=(\omega_{1}\dots,\omega_{s})$ such that, for each $j\in \{1,\dots,s\}$, $\omega_{j}\in \Omega^{p}(\cA_{j};N_{j})$ is a smooth form on $\cA_{j}$ and, for each $i,j\in \{1,\dots,s\}$, one has
$$\omega_{j_{/N_{i}\cap N_{j}}}=\omega_{i_{/N_{i}\cap N_{j}}}$$
\end{defn}






The set of all piecewise smooth forms of degree $p$ on $\underline{\cA}$ will be denoted by $\Omega^{p}(\underline{\cA};\underline{K})$. This set is a real vector space. A wedge product and an exterior derivative can be defined on $\Omega^{\ast}(\underline{\cA};\underline{K})=\bigoplus_{p\geq 0}\Omega^{p}(\underline{\cA};\underline{K})$ by the corresponding operations on each algebra $\Omega^{\ast}(\cA_{j};N_{j})=\bigoplus_{p\geq 0}\Omega^{p}(\cA_{j};N_{j})$, giving to  $\Omega^{\ast}(\underline{\cA};\underline{K})$ a structure of cochain algebra defined over $\RR$. The cohomology of this cochain algebra will be denoted by $H^{\ast}(\underline{\cA};\underline{K})$.

\vspace{3mm}

We notice that piecewise smooth cohomology of a complex of Lie algebroids defined on a simplicial complex is a particular case of this generalization. Let us briefly look at another cases made by families of manifolds on which we can define piecewise smooth cohomology. A first example is take a simplicial complex and to fix our attention on an open star of one its vertex. The family of submanifolds made by those simplices without the faces opposite to the vertex satisfies the conditions required in our definition of complex of Lie algebroids given at the beginning of this section. Any transitive Lie algebroid over the open star gives, by restriction, a complex of Lie algebroids. Another illustrative example consists of taking the family defined by intersections of open stars with any open subset of the geometric realization of a simplicial complex. The first example is obviously a particular example of this second case. The construction of a complex of Lie algebroids can be done in similar way. Our third example extends the second one and and consists of taking intersections of generalized stars with open subsets of the geometric realization of a simplicial complex. This third example is not quite different of previous examples. Nevertheless, it enhances the construction of the sheaf of the piecewise smooth forms on a complex of Lie algebroids. We provide below a description of this third example as well of the corresponding sheaf of piecewise smooth forms. Definitions and main properties of regular open subsets can be seen in \cite{maun-atop} or \cite{eile-stin}. The idea of construction of the sheaf of the piecewise smooth forms on a complex of Lie algebroids comes from \cite{penna-dgss}.


\vspace{3mm}

Let $K$ be a simplicial complex and $|K|$ its geometric realization. Consider a point $a$ of $|K|$. We recall that the generalized star of $a$, denoted also by $\St a$, is the union of the interiors of all simplices of $K$ such that $a$ belongs to those simplices (cf. with the definition 2.4.2 of \cite{maun-atop}). When the point $a$ is a vertex of $K$, it is obvious that the generalized star of $a$ is the same as the star of $a$. We notice that, for each $a\in K$, there is a unique simplex $\Delta_{a}$ of $K$ such that the point $a$ belongs to the interior of the simplex $\Delta_{a}$ (see the proposition 2.3.6 of \cite{maun-atop}).

\vspace{3mm}

\begin{prop} Let $K$ be a simplicial complex and consider a point $a$ of $|K|$. Denote by $\Delta_{a}$ the unique simplex of $K$ such that $a$ belongs to the interior of $\Delta_{a}$. Then, the generalized star of $a$ coincide with the star $\St \Delta_{a}$. Consequently, the generalized star of $a$ is an open subset of $|K|$.
\end{prop}

\proofwp. If $a$ is one of the vertices of $K$, then $\Delta_{a}=\{a\}$ and the result is proved. Suppose now that $a$ is different of any vertex of $K$. Then $a$ belongs to the interior of $\Delta_{a}$. We shall see first that $\St a \subset \St \Delta_{a}$. Let $\Delta$ be a simplex of $K$ such that $a\in \Delta$. Since $a$ is different of any vertex of $K$, it follows that $a\in \stackrel{\circ}{s}$, for some face $s$ of $\Delta$. But $a\in \ \stackrel{\circ}{\Delta_{a}}$ and so $\stackrel{\circ}{s}=\stackrel{\circ}{\Delta_{a}}$. Hence $s=\Delta_{a}$ and therefore $\Delta_{a}$ is a face of $\Delta$. We conclude then $\stackrel{\circ}{\Delta}\subset \St \Delta_{a}$. Now, let $\widetilde{\Delta}$ be a simplex of $K$ such that $\Delta_{a}$ is a face of $\widetilde{\Delta}$. Then, $a\in \widetilde{\Delta}$ and so $\stackrel{\circ}{\widetilde{\Delta}}\subset \St \Delta_{a}$. The other inclusion is obvious. The second part of the proposition is immediate. {\small $\square$}



\vspace{3mm}

Consider now $a\in |K|$ and $U$ an open subset of $|K|$ with $a\in U$. In according to \cite{penna-dgss}, the open subset $U$ is called regular open neighborhood of $a$ if $U$ is the intersection of an open neighborhood of $a$ in $|K|$ with the generalized star of the point $a$ (cf. \cite{eile-stin}, ninth section of the second chapter). Given any open subset $V$ of $|K|$, $V$ is called a regular open subset of $|K|$, if there exists a point $a\in |K|$ such that $V$ is a regular open neighborhood of the point $a$.

Obviously, a star of any simplex of a simplicial complex is a regular open subset of its geometric realization.








\vspace{3mm}

We describe now a special construction of a complex of Lie algebroids based in regular open subsets.

\vspace{3mm}

\textit{Derived complex corresponding to regular open subsets}. Let $K$ be a simplicial complex and $\underline{\cA}=\{\cA_{\Delta}\}_{\Delta\in K}$ a complex of Lie algebroids on $K$. Let $U$ be a regular open subset of $|K|$ and consider $a\in |K|$ such that $U=Z\cap \St a$, in which $Z$ is an open neighborhood of $a$ in $|K|$. Consider the unique simplex $\Delta_{a}$ of $K$ such that $a$ belongs to the interior of $\Delta_{a}$. For each simplex $\Delta\in K$ such that $\Delta_{a}$ is a face of $\Delta$, denote by $\Delta_{U}$ the set $\Delta_{U}=U\cap \Delta$.

\begin{prop} Keeping the same hypotheses and notations as above, the collection $\underline{K}^{U}$, made by the manifolds $\Delta_{U}=U\cap \Delta$ such that $\Delta_{a}$ is a face of $\Delta$, is well defined and is a transverse family.
\end{prop}

\proofwp. Let us check that the triangulation obtained in $U$ does not depend on the point $a$ chosen, that is, if  $Z\cap \St a = \widetilde{Z}\cap \St b$, then $\St a = \St b$. To see this, denote by $\Delta_{a}$ and $\Delta_{b}$ the unique simplices of $K$ which contain $a$ and $b$ in its interior respectively. Then, $\St a =\St \Delta_{a}$ and $\St b =\St \Delta_{b}$. Since $b\in V\cap \St a$, there exists a simplex $\Delta'\in K$ such that $\Delta_{a}$ is a face of $\Delta'$ and $b$ belongs to the interior of $\Delta'$. Hence, $\Delta_{b}=\Delta'$ by uniqueness of $\Delta_{b}$, and so $\Delta_{a}$ is a face of $\Delta_{b}$. Analogously, we conclude can that $\Delta_{b}$ is a face of $\Delta_{a}$ and so it holds that $\Delta_{b}=\Delta_{a}$. This show that the set $\underline{K}^{U}$ is well defined. The set $\Delta_{U}$ is a submanifold of $\Delta$. If $\Delta'$ is other simplex of $K$ such that $\Delta_{a}$ is a face of $\Delta'$, then $\Delta_{a}$ is a face of $\Delta \cap \Delta'$ and the intersection $U\cap (\Delta\cap \Delta')$ is a submanifold of $U$. This shows that the family $\underline{K}^{U}$ is transverse. {\small $\square$}

\vspace{3mm}

Keeping the same hypotheses and notations as above, the Lie algebroid $\cA_{\Delta}$ is transitive and so we can take the Lie algebroid restriction $(\cA_{\Delta})^{!!}_{\Delta_{U}}$ to $\Delta_{U}$. Therefore, we can consider the family $$\frak A_{U}=\{(\cA_{\Delta})^{!!}_{\Delta_{U}}: \Delta\in K, \ \ \Delta_{a}\prec {\Delta}\}$$ Fix two simplices $\Delta$ and $\Delta'$ of $K$ such that $\Delta_{a}$ is a common face of $\Delta$ and $\Delta'$. Let $s=\Delta\cap \Delta'$. We have $$\big((\cA_{\Delta})^{!!}_{\Delta_{U}}\big)^{!!}_{U\cap s}=(\cA_{\Delta})^{!!}_{U\cap s}=\big((\cA_{\Delta})^{!!}_{s}\big)^{!!}_{U\cap s}=(\cA_{s})^{!!}_{U\cap s}$$ and analogously $$\big((\cA_{\Delta'})^{!!}_{\Delta'_{U}}\big)^{!!}_{U\cap s}=(\cA_{s})^{!!}_{U\cap s}$$ Hence, we proved the following proposition.

\vspace{3mm}

\begin{prop} The family $\frak A_{U}$ is a complex of Lie algebroids defined over the family of manifolds $\underline{K}^{U}$. The complex $\frak A_{U}$ is called the derived complex of the complex $\underline{\cA}$ corresponding to the regular open subset $U$.
\end{prop}


The cochain algebra of the piecewise smooth forms on the derived complex of a complex of Lie algebroids will be denoted simply by $\Omega^{\ast}_{ps}(\frak A_{U})$.

\vspace{3mm}

Keeping the same hypotheses and notations as above, let $U$ and $V$ two regular open subsets of $|K|$ such that $V\subset U$. We shall construct now a restriction mapping from $\Omega^{\ast}_{ps}(\frak A_{U})$ to $\Omega^{\ast}_{ps}(\frak A_{V})$. Let $a$ and $b$ $\in |K|$ such that $U$ and $V$ are regular open neighborhoods of $a$ and $b$ respectively. Denote by $\Delta_{a}$ and $\Delta_{b}$ the unique simplices of $K$ which contain $a$ and $b$ in its interior respectively. Since $b\in U$, there exists a simplex $\Delta \in K$ such that $\Delta_{a}$ is a face of $\Delta$ and the point $b$ belongs to the interior of $\Delta$. Hence, $\Delta_{b}=\Delta$ and so $\Delta_{a}$ is a face of $\Delta_{b}$. If $\Delta'$ is a simplex of $K$ such that $\Delta_{b}$ is a face $\Delta'$, then $\Delta_{a}$ is a face of $\Delta'$ and, consequently, every $V\cap \Delta$ of $\underline{K}^{V}=\{\Delta_{V}: \Delta\in K, \Delta_{b}\prec \Delta\}$ belongs to $\underline{K}^{U}=\{\Delta_{U}:\Delta\in K, \Delta_{a}\prec \Delta\}$. Let $\omega=(\omega_{\Delta_{U}})_{\Delta_{U}\in \underline{K}^{U}}$ $\in \Omega^{\ast}_{ps}(\frak A_{U})$ be a piecewise sooth form on the complex of Lie algebroids $\frak A_{U}$. For each simplex $\Delta\in K$ such that $\Delta_{b}$ is a face of $\Delta$, we have that $\big((\cA_{\Delta})^{!!}_{\Delta_{U}}\big)^{!!}_{\Delta_{V}}=(\cA_{\Delta})^{!!}_{\Delta_{V}}$ and we can restrict the smooth form $$\omega_{\Delta_{U}}\in \Omega^{\ast}\big((\cA_{\Delta})^{!!}_{\Delta_{U}};\Delta_{U}\big)$$ to the submanifold $\Delta_{V}$, obtaining the smooth form $$\omega_{\Delta_{V}}=(\omega_{\Delta_{U}})^{!!}_{\Delta_{V}}\in \Omega^{\ast}\big((\cA_{\Delta})^{!!}_{\Delta_{V}};\Delta_{V}\big)$$ Therefore, we obtain the differential form $(\omega_{\Delta_{V}})_{\Delta_{V}\in \underline{K}^{V}}$. Basic arguments show that the form $(\omega_{\Delta_{V}})_{\Delta_{V}\in \underline{K}^{V}}$ is a piecewise smooth form and so it belongs to $\in \Omega^{\ast}_{ps}(\frak A_{V})$. As done before, the piecewise smooth form $(\omega_{\Delta_{V}})_{\Delta_{V}\in \underline{K}^{V}}$ is denoted by $\omega_{/V}$.

\vspace{3mm}

The following proposition is easily derived.

\vspace{3mm}

\begin{thm} Let $K$ be a simplicial complex and $\underline{\cA}=\{\cA_{\Delta}\}_{\Delta\in K}$ a complex of Lie algebroids on $K$. Let $U$ and $V$ be two regular open subsets of $|K|$ such that $U\subset V$ and consider the derived complexes of Lie algebroids $\frak A_{U}$ and $\frak A_{V}$ corresponding to $U$ and $V$ respectively. For each $p\geq 0$, denote by $$r^{p^{K}}_{_{L}}:\Omega^{p}_{ps}(\frak A_{U})\longrightarrow \Omega^{p}_{ps}(\frak A_{V})$$ the map induced by restriction, that is, for each $\omega\in \Omega^{p}_{ps}(\frak A_{U})$,  $$r^{p^{U}}_{_{V}}(\omega)=\omega_{/V}$$

\begin{itemize}

\item For each $p\geq 0$,  $r^{p^{U}}_{_{U}}=id_{\Omega^{p}_{ps}(\frak A_{U})}$.

\item $r^{\ast^{U}}_{_{V}}:\Omega^{\ast}_{ps}(\frak A_{U})\longrightarrow \Omega^{\ast}_{ps}(\frak A_{V})$ is a morphism of graded algebras.

\item If $W$ is other generalized regular open subset of $|K|$ with $W\subset V$ and $\frak A_{W}$ is the derived complex of Lie algebroids corresponding
to $W$, then the diagram below is a commutative diagram of cochain complexes

$$
\xymatrix{
\Omega^{\ast}_{ps}(\frak A_{U})\ar[dr]_{r^{\ast^{U}}_{_{W}}}\ar[rr]^{r^{\ast^{U}}_{_{V}}} & & \Omega^{\ast}_{ps}(\frak A_{V})\ar[dl]^{r^{\ast^{V}}_{_{W}}} \\
& \Omega^{\ast}_{ps}(\frak A_{W}) &
}
$$

\end{itemize}
\end{thm}

\vspace{3mm}

A direct consequence from the previous proposition is the following corollary.

\vspace{3mm} 

\begin{cor} Keeping the same hypotheses and notations from the previous theorem, for each $p\geq 0$, the correspondence which associates, to each regular open subset $U$ of $|K|$, the real vector space $\Omega^{p}_{ps}(\frak A_{U})$ of the piecewise smooth forms defined on $U$ and, to each pair of regular open subsets $U$ and $V$ of $|K|$ such that $V\subset U$, the homomorphism $r^{p^{U}}_{_{V}}$, is a presheaf. This presheaf is called the presheaf of the piecewise smooth forms of degree $p$ on the complex of Lie algebroids $\underline{\cA}$.
\end{cor}

\vspace{3mm}

The previous proposition leads us to the following definition (cf. with \cite{penna-dgss}).

\vspace{3mm}

\begin{defn} Let $K$ be a simplicial complex and $\underline{\cA}=\{\cA_{\Delta}\}_{\Delta\in K}$ a complex of Lie algebroids on $K$. For each $p\geq 0$, the sheaf of the piecewise smooth forms of degree $p$ on the complex of Lie algebroids $\underline{\cA}$ is the sheaf constructed canonically from the presheaf of the piecewise smooth forms of degree $p$ on the complex of Lie algebroids $\underline{\cA}$.
\end{defn}

\vspace{3mm}

We state now the last result of this paper.

\vspace{3mm}

\begin{thm} Let $K$ be a simplicial complex and $\underline{\cA}=\{\cA_{\Delta}\}_{\Delta\in K}$ a complex of Lie algebroids on $K$. Then the sheaf $\frak S$ of the piecewise smooth forms of degree $p$ on the complex of Lie algebroids $\underline{\cA}$ is fine.
\end{thm}

\proofwp. Let $\frak U=\{U_{j}\}_{j\in J}$ be a locally finite open covering of $|K|$ by regular open subsets of $|K|$. Since the set of all regular open subsets of $|K|$ is a base for the topology of $|K|$, we can assume that each open subset $U\in \frak U$ is a regular open subset. If $\{\varphi_{j}\}_{j\in J}$ is piecewise smooth partition of unity subordinated to the covering $\frak U$, the homomorphisms of presheaves $h_{j}: \Omega^{p}_{ps}(\frak A_{U}) \longrightarrow \Omega^{p}_{ps}(\frak A_{U})$ defined by $h_{j}(\omega)=\varphi_{j_{/U}}\omega$ for each $\omega\in \Omega^{p}_{ps}(\frak A_{U})$ induce homomorphisms from $\frak A^{p}$ to $\frak A^{p}$ satisfying the conditions  which characterize the definition of fine sheaf. Therefore, the result is proved if we find a piecewise smooth partition of unity subordinated to the covering $\frak U$. By lemma shrinking, there is an open covering $\frak V=\{V_{j}\}_{j\in J}$ such that, for each $j\in J$, $\overline{V_{j}}\subset U_{j}$. Let $U\in \frak U$ and $V\in \frak V$ such that  $\overline{V}\subset U$. Consider $a\in |K|$ such that $U$ is a regular open neighborhood of $a$ in $|K|$. For each simplex $\Delta\in K$ such that $\Delta_{a}$ is a face of $\Delta$, consider the closed subset $\overline{V}\cap \Delta$ of $\Delta$. Take the union of all $\overline{V}\cap \Delta$ such that $\Delta_{a}$ is a face of $\Delta$ and denote that union by $W$. Since $|K|$ is compact, the topology of $|K|$ coincide with the topology induced from the Euclidian space. We have that $W$ a closed subset of the Euclidian space. The open star $\St \Delta_{a}$ is open in $|K|$ and so there is an open subset $Z$ of the Euclidian space such that $\St \Delta_{a}=Z\cap |K|$. The closed subset is contained in the open subset $Z$. Hence, we can fix a smooth function $\varphi:Z\longrightarrow \RR$ such that $\varphi$ does not vanish on $W$. By restriction to each submanifold $\Delta_{U}=U\cap \Delta$, we have a piecewise smooth function on $U$ which does not vanish on each $\overline{V}\cap \Delta$. Take the sum of these functions and consider the quotient of each function by the sum. This defines a partition of unity made by piecewise smooth functions. {\small $\square$}

\vspace{3mm}


\end{document}